\documentclass[12pt,a4paper]{article}
\usepackage{amsmath,amsfonts,amssymb,amscd}

\usepackage{graphicx}
\usepackage{amsfonts,amssymb,amscd,amsmath}
\usepackage[body={6.6in, 9.2in},left=1in, top=1.2in]{geometry}
\usepackage{indentfirst}

\usepackage{tikz}
\usetikzlibrary{matrix}

\def\H{\operatorname{H}}

\def\Ker{\operatorname{Ker}}
\def\Ad{\operatorname{Ad}}

\def\ph{\varphi}

\def\d{\operatorname{d}}
\def\GL{\operatorname{GL}}

\def\Q{\operatorname{Q}}

\def\Lie{\operatorname{Lie}}

\newcounter{th}
\def\t{\refstepcounter{th}{\bf \noindent{Theorem} \arabic{th}. }}

\newcounter{prop}
\def\prop{\refstepcounter{prop}{\bf \noindent{Proposition} \arabic{prop}. }}

\newcounter{lem}
\def\lem{\refstepcounter{lem}{\bf \noindent{Lemma} \arabic{lem}. }}

\newcounter{de}

\newcounter{ex}

\begin{document}

\begin{center}

{\Large{\bf Vector fields on $\Pi$-symmetric flag supermanifolds}}\footnote[1]{Supported by Max Planck Institute for Mathematics, Bonn.}
\bigskip

{\bf E.G.Vishnyakova}\\[0.3cm]
\end{center}

\bigskip

\begin{abstract}
 The main result of this paper is the computation of the Lie superalgebras of holomorphic
vector fields on the complex $\Pi$-symmetric flag supermanifolds, introduced by Yu.I.~Manin. We prove  that with one exception any vector field is fundamental with respect to the
natural action of the Lie superalgebra $\mathfrak q_n(\mathbb C)$.

\end{abstract}

\bigskip

\section{Introduction}

A $\Pi$-symmetric flag super\-manifold is a subsupermanifold in a flag supermanifold in $\mathbb C^{n|n}$ that is invariant with respect to an odd involution in $\mathbb C^{n|n}$. This supermanifold possesses a transitive action of the linear classical Lie superalgebra $\mathfrak q_n(\mathbb C)$,  which belongs to one of  two ``strange series'' in the Kac classification \cite{Kac}. It turns out that with one exceptional case all global holomorphic vector fields are fundamental for this action of the Lie superalgebra $\mathfrak q_n(\mathbb C)$. In the simplest case of super-Grassmannians the similar result  was obtained in \cite{Oni}. 

The main result of this paper was announced in \cite{Vivector} and the idea of the proof was given in \cite{ViPi-sym}.  The goal of this notes is to give a detailed proof. 
We also describe the connected component of the automorphism supergroup of this supermanifolds.

\section{Flag supermanifolds}

We will use the word ``supermanifold'' \,in the sense of Berezin and Leites \cite{BL}, see also \cite{Oni} for details. Throughout, we will restrict our attention to the complex-analytic version of the theory of supermanifolds. Recall that a {\it complex-analytic superdomain of dimension $n|m$} is a $\mathbb{Z}_2$-graded
ringed space of the form $(\mathcal U_0, \mathcal{F}_{\mathcal U_0} \otimes_{\mathbb C} \bigwedge(m) )$,
where $\mathcal{F}_{\mathcal U_0}$ is the sheaf of holomorphic functions on an open set $\mathcal U_0\subset \mathbb{C}^n$ and
$ \bigwedge(m)$ is the exterior (or Grassmann) algebra with $m$ generators.
 A {\it complex-analytic supermanifold} of dimension $n|m$ is a $\mathbb{Z}_2$-graded locally ringed space that
is locally isomorphic to a complex superdomain of dimension $n|m$. Let $\mathcal{M} = (\mathcal{M}_0,{\mathcal O}_{\mathcal{M}})$ be a supermanifold and $\mathcal{J}$ be the subsheaf of ideals generated by odd elements in ${\mathcal O}_{\mathcal{M}}$. We set $\mathcal{F}_{\mathcal{M}_0}:= {\mathcal
	O}_{\mathcal{M}}/\mathcal{J}$. Then $(\mathcal{M}_0,
\mathcal{F}_{\mathcal{M}_0})$ is a usual complex-analytic manifold, it
is called the {\it underlying space} of $\mathcal{M}$. Usually we will
write $\mathcal{M}_0$ instead of $(\mathcal{M}_0,
\mathcal{F}_{\mathcal{M}_0})$.

In this paper we denote by $\mathbf
F_{k|l}^{m|n}$ a flag supermanifold of type $k|l$ in the vector superspace $\mathbb C^{m|n}$. Here we set $k=(k_1,\ldots,k_r)$ and
$l=(l_1,\ldots,l_r)$ such that
\begin{equation}\label{eq conditions on k_i and l_j}
\begin{split}
0\le k_r\le\ldots\le k_1\le m&,\quad 0\le l_r\ldots\le l_1\le n\quad \text{and}\\
 \quad 0 <
k_r +l_r <&\ldots < k_1 + l_1 < m+n.
\end{split}
\end{equation} 
 A $\Pi$-symmetric flag supermanifold $\mathbf{\Pi F}_{k|k}^{n|n}$ of type $k=(k_1,\ldots,k_r)$ in $\mathbb C^{n|n}$ is a certain subsuper\-manifold in $\mathbf F_{k|k}^{n|n}$. Let us give an explicite description of these supermanifolds in terms of charts and local coordinates (see also \cite{Man,Vi,Viholom}). 

Let us take two non-negative integers $m,n\in\mathbb Z$ and two sets of non-negative integers $$
k=(k_1,\ldots,k_r),\quad  \text{and} \quad  l=(l_1,\ldots,l_r)
$$
such that (\ref{eq conditions on k_i and l_j}) holds. 
The underlying space of the supermanifold
$\mathbf F_{k|l}^{m|n}$ is the product  $\mathbf
F^{m}_{k}\times \mathbf F^{n}_{l}$ of two manifolds of flags of type $k=(k_1,\ldots,k_r)$ and $l=(l_1,\ldots,l_r)$  in $\mathbb C^m$ and $\mathbb C^n$, respectively. For any $s = 1,\ldots,r$ let us fix two subsets 
$$
I_{s\bar
0}\subset\{1,\ldots,k_{s-1}\}\quad \text{and} \quad I_{s\bar 1}\subset\{1,\ldots,l_{s-1}\},
$$
 where $k_0=m$ and $l_0=n$, such that
$|I_{s\bar 0}| = k_s,$ and $|I_{s\bar 1}| = l_s$. We set
$I_s=(I_{s\bar 0},I_{s\bar 1})$ and $I = (I_1,\ldots,I_r)$. 
Let us assign the following $(k_{s-1} + l_{s-1})\times (k_s + l_s)$-matrix 
\begin{equation}\label{eq local chart general}
Z_{I_s}=\left(
\begin{array}{cc}
X_s & \Xi_s\\
\H_s & Y_s \end{array} \right), \ \ s=1,\dots,r,
\end{equation}
to any $I_s$. Here we assume that $X_s=(x^s_{ij})\in \operatorname{Mat}_{k_{s-1}\times k_s}(\mathbb C),\;
Y_s=(y^s_{ij})\in \operatorname{Mat}_{l_{s-1}\times l_s}(\mathbb C)$, and elements of the matrices 
$\Xi_s=(\xi^s_{ij})$, $\H_s=(\eta^s_{ij})$ are odd. We also assume that $Z_{I_s}$ contains the identity submatrix $E_{k_s+l_s}$ of size $(k_s+l_s)\times (k_s+l_s)$ in the lines with numbers $i\in I_{s\bar 0}$ and $k_{s-1} + i,\; i\in I_{s\bar 1}$. For example in case 
$$
I_{s\bar 0}=\{k_{s-1}-k_s+1,\ldots, k_{s-1}\}\quad  \text{and} \quad I_{s\bar 1}=\{l_{s-1}-l_s+1,\ldots, l_{s-1}\}
$$
the matrix $Z_{I_s}$ has the following form:
$$
Z_{I_1} =\left(
\begin{array}{cc}
X_s&\Xi_s\\
E_{k_s}&0\\
\H_s&Y_s\\0&E_{k_s}\end{array} \right).
$$
(For simplisity of notation we use here the same letters $X_s$, $Y_s$, $\Xi_s$ and $\H_s$ as in (\ref{eq local chart general}).)

We see that the sets $I_{\bar 0} =(I_{1\bar 0},\ldots,I_{r\bar 0})$ and $I_{\bar
1}= (I_{1\bar 1},\ldots,I_{r\bar 1})$ determine the charts $U_{I_{\bar 0}}$ and $V_{I_{\bar 1}}$ on the flag manifolds $\mathbf F^{m}_{k}$ and $\mathbf F^{n}_{l}$, respectively. We can take the non-trivial elements (i.e., those not contained in the identity submatrix) from
 $X_s$ and $Y_s$ as local coordinates in $U_{I_{\bar 0}}$ and $U_{I_{\bar 1}}$, respectively. Summing up, we defined an atlas 
$$
\{U_I = U_{I_{\bar 0}}\times
U_{I_{\bar 1}}\} \quad \text{on} \quad \mathbf F^{m}_{k} \times\mathbf
F^{n}_{l}
$$
 with chards parametrized by $I = (I_s)$. In addition the sets $I_{\bar 0}$ and $I_{\bar 1}$ determine the superdomain $\mathcal U_I$ with underlying space $U_I$ and with even and odd coordinates $x^s_{ij}$, $y^s_{ij}$ and $\xi^s_{ij}$, $\eta^s_{ij}$, respectively. (As above we assume that $x^s_{ij}$, $y^s_{ij}$, $\xi^s_{ij}$ and $\eta^s_{ij}$ are non-trivial. That is they are  not contained in the identity submatrix.)
Let us define the transition functions between two  superdomains corresponding to
 $I = (I_s)$ and $J = (J_s)$ by the following formulas:
\begin{equation}\label{eq transition functions}
Z_{J_1} = Z_{I_1}C_{I_1J_1}^{-1}, \quad Z_{J_s} =
C_{I_{s-1}J_{s-1}}Z_{I_s}C_{I_sJ_s}^{-1},\quad s\ge 2. 
\end{equation}
Here $C_{I_sJ_s}$ is an invertible submatrix in $Z_{I_s}$ that coinsists of 
lines with numbers $i\in J_{s\bar 0}$ and $k_{s-1} + i,$ where $i\in J_{s\bar 1}$. In other words, we choose the matrix $C_{I_sJ_s}$ in such a way that $Z_{J_s}$ contains the identity submatrix $E_{k_s+l_s}$ in lines with numbers $i\in J_{s\bar 0}$ and $k_{s-1} + i,$ where $i\in J_{s\bar 1}$.
These charts and transition functions define a supermanifold that we denote by
$\mathbf F_{k|l}^{m|n}$. This supermanifold we will call the {\it supermanifold of flags} of type $k|l$. In case $r = 1$ this supermanifold is called the {\it super-Grass\-mannian} and is denoted by $\mathbf
{Gr}_{m|n,k|l}$ (see also \cite{Oni,Man}).

Let us take $n\in\mathbb N$ and $k=(k_1,\ldots,k_r)$, such that 
$$
0 < k_1 <\ldots < k_r < n.
$$
We will define the {\it supermanifold of $\Pi$-symmetric flags $\mathbf{\Pi F}_{k|k}^{n|n}$ of type}
$k$ in $\mathbb C^{n|n}$ as a certain subsupermanifold in $\mathbf F_{k|k}^{n|n}$. The underlying space of $\mathbf{\Pi F}_{k|k}^{n|n}$ is the diagonal in $\mathbf
F^{n}_{k}\times\mathbf F^{n}_{k}$, that is clearly isomorphic to $\mathbf F^{n}_{k}$. For any $s = 1,\ldots,r$ we fix a set $I_{s\bar 0}= I_{s\bar 1}\subset\{1,\ldots,k_{s-1}\}$,
where $|I_{s\bar 0}| = k_s$ and $k_0=n$. Consider the chart on $\mathbf F_{k|k}^{n|n}$ corresponding to $I= (I_s)$, where $I_s = (I_{s\bar 0},I_{s\bar 0})$. Such charts cover the diagonal in $\mathbf F^{n}_{k}\times\mathbf F^{n}_{k}$. Let us define the subsupermanifold of $\Pi$-symmetric flags in these charts by the equations $X_s=Y_s$, $\Xi_s=\H_s$. It is easy to see that these equations are well-defined with respect to the transition functions (\ref{eq transition functions}). The coordinate matrices in this case have the following form 
\begin{equation}\label{eq coordinats matrices P-sym}
Z_{I_s}=\left(
\begin{array}{cc}
X_s & \Xi_s\\
\Xi_s & X_s \end{array} \right), \ \ s=1,\dots,r.
\end{equation}
(Compare with  (\ref{eq local chart general}).) As above even and odd local coordinates on $\mathbf{\Pi F}_{k|k}^{n|n}$ are non-trivial elements from $X_s$ and  $\Xi_s$, respectively.
The transition functions between two charts are defined again by formulas (\ref{eq transition functions}). We can consider the supermanifold $\mathbf{\Pi F}_{k|k}^{n|n}$ as the ``set of 	fixed-point'' of a certain odd involution ${\Pi}$ in $\mathbb C^{n|n}$ (see
\cite{Man}). In case $r = 1$, the supermanifold of  $\Pi$-symmetric flags is called also the
 {\it $\Pi$-symmetric super-Grassmannian}. We will denote it by $\mathbf {\Pi Gr}_{n|n,k|k}$.

Let $\mathcal M = (\mathcal M_0,\mathcal O_{\mathcal M})$ be a complex-analytic supermanifold. 
Denote by 
$\mathcal T = \mathcal Der\,(\mathcal O_{\mathcal M})$ the tangent sheaf or the sheaf of vector fields on  $\mathcal M$. It is a sheaf of Lie superalgebras with respect to the multiplication $[X,Y] = YX
-(-1)^{p(X)p(Y)}XY$. The global sections of $\mathcal T$
are called {\it holomorphic vector fields} on $\mathcal M$. They form a complex Lie superalgebra
that we will denote by $\mathfrak v(\mathcal M)$. This Lie superalgebra is finite dimensional if $\mathcal M_0$ is compact. The goal of this paper is to compute the Lie superalgebra $\mathfrak v(\mathcal M)$ in the case when $\mathcal M$ is a supermanifold of $\Pi$-symmetric flags of type $k$ in $\mathbb C^{n|n}$.

We denote by $\mathfrak q_n(\mathbb C)$ the Lie subsuperalgebra in $\mathfrak {gl}_{n|n}(\mathbb C)$ that coinsists of the following marices:
$$
\left(
  \begin{array}{cc}
    A & B \\
    B & A \\
  \end{array}
\right), \quad\text{where} \quad  A,B\in \mathfrak{gl}_n(\mathbb C).
$$
Denote by $\Q_{n}(\mathbb C)$ the Lie supergroup of $\mathfrak q_n(\mathbb C)$. In \cite{Man} an action of  $\Q_{n}(\mathbb C)$ on the supermanifold
 $\mathbf{\Pi F}_{k|k}^{n|n}$ was defined. In our coordinates this action is given by the following formulas:
\begin{equation}\label{eq action of Q}
\begin{aligned} 
&(L,(Z_{I_1},\ldots,Z_{I_r})) \longmapsto
(\tilde Z_{J_1},\ldots,\tilde Z_{J_r}), \quad \text{where}\\
\ \ &L\in \Q_{n}(\mathbb C),\quad \tilde Z_{J_1} =
LZ_{I_1}C_1^{-1},\quad \tilde Z_{J_s} = C_{s-1}Z_{I_s}C_s^{-1}.
\end{aligned} 
\end{equation}
Here $C_1$ is an invertible submatrix in $LZ_{I_1}$ that consists of lines with numbers $i$ and $n+i$, where  $i\in J_1$; and $C_s,\; s\ge 2$, is an invertible submatrix in $C_{s-1}Z_{I_s}$ that consists of lines with numbers $i$ and $k_{s-1}+i$, where
$i\in J_s$. This Lie supergroup action induces a Lie superalgebra homomorphism 
$$
\mu:\mathfrak{q}_{n}(\mathbb C)\to\mathfrak v(\mathbf {\Pi F}_{k|k}^{n|n}).
$$ 
In case
$r=1$ in \cite[Proposition 5.5]{Oni} it was proven that $\Ker\mu = \langle E_{2n}\rangle$, where
$E_{2n}$ is the identity matrix of size $2n$. In general case $r>1$ the proof is similar.  Hence, $\mu$ induces an injective homomorphism of Lie superalgebras $\mathfrak{q}_{n}(\mathbb
C)/\langle E_{2n}\rangle\to\mathfrak v(\mathbf {\Pi F}_{k|k}^{n|n})$. We will show that with one exception this homomorphism is an isomorphism.

\section{About superbundles}

Recall that a {\it morphism} of a complex-analytic supermanifold 
$\mathcal M$ to a complex-analytic supermanifold $\mathcal N$ is a pair $f = (f_0,f^*)$, where $f_0: \mathcal M_0\to \mathcal N_0$ is a holomorphic map and $f^*: \mathcal O_{\mathcal N}\to (f_0)_*(\mathcal
O_{\mathcal M})$ is a homomorphism of sheaves of superalgebras. 

\medskip
\noindent {\bf Definition.} We say that a {\it
superbundle} with fiber $\mathcal S$, base $\mathcal
B$, total space $\mathcal M$ and projection $p=(p_0,p^*):
\mathcal M\rightarrow \mathcal B$ is given if 
there exists an open covering 
 $\{U_i\}$ on $\mathcal B_0$ and isomorphisms $\psi_i: (p_0^{-1}(U_i),\mathcal O_{\mathcal M})\rightarrow
(U_i,\mathcal {\mathcal O}_{\mathcal B})\times \mathcal S$ such that the following diagram is commutative:
\begin{center}
\begin{tikzpicture}
\matrix (m) [matrix of math nodes,row sep=3em,column sep=4em,minimum width=2em]
 {(p_0^{-1}(U_i),\mathcal O_{\mathcal M}) & (U_i,\mathcal {\mathcal O}_{\mathcal B})\times \mathcal S \\
	(U_i,\mathcal {\mathcal O}_{\mathcal B}) & (U_i,\mathcal {\mathcal O}_{\mathcal B}) \\};
\path[-stealth]
(m-1-1) edge node [left] {$p$} (m-2-1)
edge node [below] {$\psi_i$} (m-1-2)
(m-2-1.east|-m-2-2) edge node [below] {$id$}  (m-2-2)
(m-1-2) edge node [right] {$pr$} (m-2-2);
\end{tikzpicture}
\end{center}
where $pr$ is the natural projection.
\medskip

\noindent {\bf Remark.} From the form of transition functions (\ref{eq transition functions}) it follows that for $r
> 1$ the supermanifold $\mathbf{\Pi F}_{k|k}^{n|n}$ is a superbundle with base $\mathbf {\Pi Gr}_{n|n,k_1|k_1}$
and fiber $\mathbf {\Pi F}_{k'|k'}^{k_1|k_1}$, where $k'=(k_2,\ldots,k_r)$. In local coordinates the projection $p$ is given by 
$$
(Z_1,Z_2,\ldots Z_n) \longmapsto (Z_1).
$$
Moreover, the formulas (\ref{eq action of Q}) tell us that the projection $p$ is equivariant with respect to the action of the supergroup  $\Q_{n}(\mathbb C)$ on $\mathbf{\Pi F}_{k|k}^{n|n}$ and $\mathbf {\Pi Gr}_{n|n,k_1|k_1}$.
\medskip

Let $p=(p_0,p^*): \mathcal M\to \mathcal N$ be a morphism of supermanifolds. A vector field $v\in\mathfrak
v(\mathcal M)$ is called {\it projectable} with respect to $p$,
if there exists a vector field $v_1\in\mathfrak v(\mathcal N)$
such that 
$$
p^*(v_1(f))=v(p^*(f))\quad \text{for all} \quad f\in \mathcal O_{\mathcal N}.
$$
In this case we say that {\it $v$ is projected into} $v_1$.
Projectable vector fields form a Lie subsuperalgebra  $\overline{\mathfrak
v}(\mathcal M)$ in $\mathfrak v(\mathcal M)$.
In case if $p$ is a projection of a superbundle, the homomorphism
$p^*: \mathcal O_{\mathcal N}\to p_*(\mathcal O_{\mathcal M})$ is injective. Hence, any projectable vector field $v$ is projected into unique vector field
 $v_1 = \mathcal{P}(v)$. The map 
 $$
 \mathcal{P}:
\overline{\mathfrak v}(\mathcal M)\to \mathfrak v(\mathcal N),\quad v \mapsto v_1
$$
is a homomorphism of Lie superalgebras. A vector field $v\in\mathfrak v(\mathcal M)$ is called {\it vertical}, if $\mathcal{P}(v)=0$.
Vertical vector fields form an ideal  $\Ker\mathcal{P}$ in
Ëè $\overline{\mathfrak v}(\mathcal M)$.

We will need the following proposition proved in \cite{Bash}.

\medskip
\prop\label{prop Bash} {\sl Let 
	$p: \mathcal M\to
\mathcal B$ be the projection of a superbundle with fiber $\mathcal S$. If $\mathcal O_{\mathcal S}(\mathcal S_0) = \mathbb C$, then any holomorphic vector field from $\mathfrak v(\mathcal M)$ is projectable with respect to $p$.}
\medskip

For any superbundle  $p: \mathcal M\to
\mathcal B$ with fiber $\mathcal S$ we define the sheaf $\mathcal W$ on
$\mathcal B_0$ in the following way. We asign to any open set $U\subset \mathcal B_0$
 the set of all vertical vector fields on the supermanifold $(p_0^{-1}(U),\mathcal
O_{\mathcal M})$. In \cite{Vi} the following statement was proven.

\medskip
\prop\label{prop W is localy free}  {\sl Assume that $\mathcal S_0$ is compact. Then $\mathcal W$ is a localy free 
sheaf of 
$\mathcal O_{\mathcal B}$-modules and
 $\dim\mathcal
W= \dim\mathfrak v(\mathcal S)$. The Lie algebra $\mathcal W(\mathcal B_0)$
coincides with the ideal of all vertical vector fields in $\mathfrak
v(\mathcal M)$}.
\medskip

Let us describe the corresponding to $\mathcal W$ graded sheaf as in \cite{Vi}.  Consider the following filtration in ${\mathcal O}_{\mathcal B}$
$$
{\mathcal O}_{\mathcal B} = \mathcal J^0 \supset \mathcal J^1 \supset \mathcal
J^2\dots
$$ 
where $\mathcal J$ is the sheaf of ideals in ${\mathcal O}_{\mathcal B}$ generated by odd elements. We have the corresponding graded sheaf of superalgebras
$$
\tilde{\mathcal O}_{\mathcal B}=\bigoplus_{p\ge 0}(\tilde{\mathcal O}_{\mathcal B})_{p},
$$
where  $(\tilde{\mathcal O}_{\mathcal B})_{p}=\mathcal J^p/\mathcal J^{p+1}$.
Putting $\mathcal W_{(p)}=\mathcal J^p\mathcal W$ we get the following filtration in $\mathcal W$:
\begin{equation}\label{eq filtrations in W}
\mathcal W=\mathcal W_{(0)}\supset\mathcal W_{(1)}\supset\dots .
\end{equation}
We define the $\mathbb Z$-graded sheaf of $\mathcal F_{\mathcal B_0}$-modules by
\begin{equation}\label{eq tilda W def}
\tilde {\mathcal W}=\bigoplus_{p\ge 0}\tilde{\mathcal W}_{p}, \ \
\text{where}\ \ \tilde{\mathcal W}_{p}=\mathcal W_{(p)}/\mathcal
W_{(p+1)}, 
\end{equation}
where $\mathcal F_{\mathcal B_0}$ is the structure sheaf of the underlying space $\mathcal B_0$. 
The ${\mathbb Z}_2$-grading in $\mathcal W_{(p)}$ induces the
${\mathbb Z}_2$-grading in $\tilde{\mathcal W}_{p}$. Note that the natural map $\mathcal W_{(p)} \to \tilde{\mathcal
W}_{p}$ is even.

\section{Functions on $\Pi$-symmetric flag supermanifolds}

In this section we show that the superbundle described in Section 2, that is the
$\Pi$-symmetric flag supermanifold, satisfies conditions of Proposition \ref{prop Bash}. Holomorphic functions on other flag supermanifolds was considered in \cite{Viholom}. 

\medskip
\lem \label{lem constant functions on fiber and base} {\sl Let $\mathcal M$ be a superbundle with base  $\mathcal B$ and fiber
$\mathcal S$. Assume that $\mathcal O_{\mathcal B}(\mathcal B_0) = \mathbb C$ and $\mathcal
O_{\mathcal S}(\mathcal S_0) =\mathbb C$. Then $\mathcal O_{\mathcal M}(\mathcal M_0) =\mathbb C$.}
\medskip

In the Lie superalgebra  $\mathfrak {q}_{n}(\mathbb C)_{\bar 0}\simeq \mathfrak{gl}_n(\mathbb C)$
we fix the following Cartan subalgebra: 
$$
\mathfrak t=
\{\operatorname{diag}(\mu_1,\dots,\mu_n)\},
$$
 the following system of positive roots:
$$
\Delta^+=\{\mu_i-\mu_j, \,\,i<j\}
$$
and the following system of simple roots:
$$
\begin{array}{c}
\Phi= \{\alpha_1,..., \alpha_{n-1}\}, \,\,\,
\alpha_i=\mu_i-\mu_{i+1}.
\end{array}
$$

Denote by $\mathfrak t^*(\mathbb R)$
a real subspace in $\mathfrak t^*$
spaned by $\mu_j$. Consider the scalar product $( \,,\, )$ in $\mathfrak t^*(\mathbb R)$ such that the vectors $\mu_j$ form an orthonormal basis. An element $\gamma\in \mathfrak t^*(\mathbb R)$ is called {\it dominant} if $(\gamma, \alpha)\ge 0$ for all $\alpha \in \Delta^+$.

We need the Borel-Weyl-Bott Theorem (see for example \cite{ADima} for details). Let  $G\simeq \GL_{n}(\mathbb C)$ be the underlying space of $\Q_n(\mathbb C)$,  $P$ be a parabolic subgroup in $G$ and $R$ be the reductive part of $P$. Assume that $\mathbf E_{\varphi}\to G/P$ is the homogeneous vector bundle corresponding to a representation  $\varphi$ of
$P$ in $E=(\mathbf E_{\varphi})_{P}$. Denote by $\mathcal E_{\varphi}$ the sheaf of holomorphic section of this vector bundle.

\medskip

\t [Borel-Weyl-Bott]. \label{teor borel} {\sl Assume that the representation	$\varphi: P\to GL(E)$ is completely reducible and $\lambda_1,..., \lambda_s$
	are highest weights of $\varphi|R$. Then the $G$-module $H^0(G/P,\mathcal E_{\varphi})$ is isomorphic to the sum of irruducible $G$-modules with highest weights $\lambda_{i_1},..., \lambda_{i_t}$, where 
	$\lambda_{i_a}$ are dominant highest weights.

}

\medskip

The main result of this section is the following theorem.

\medskip
\t\label{teor constant functions} {\sl Let $\mathcal M=\mathbf{\Pi F}_{k|k}^{n|n}$, then $\mathcal O_{\mathcal M}(\mathcal M_0) =\mathbb C$.}
\medskip

\noindent{\it Proof.} First consider the case $r=1$. This is $\mathcal M= \mathbf{\Pi Gr}_{n|n,k|k}$. Let us prove that 
$\tilde{\mathcal O}_{\mathcal M}(\mathcal M_0)=\mathbb C$, where $\tilde{\mathcal O}_{\mathcal M}$ is defined as in the previous section. We use the Borel-Weyl-Bott Theorem. The manifold $\mathcal M_0=
\mathbf{Gr}_{n,k}$ is isomorphic to $G/H$, where
$G=\GL_{n}(\mathbb C)$ and
\begin{equation}\label{eq_stab Grassmannian}
H=  \left\{\left.\left(
\begin{array}{cc}
A & 0\\
C& B\end{array} \right)\right|\,\, A\in \GL_{n-k}(\mathbb C), \,\,\,B\in
\GL_{k}(\mathbb C) \right\}. 
\end{equation}
The reductive part $R$ of $H$ is given by the equation $C=0$. It is isomorphic to $\GL_{n-k}(\mathbb C)\times \GL_{k}(\mathbb C)$. Let
$\rho_{1}$ and $\rho_{2}$ be the standard representations of
$\GL_{n-k}(\mathbb C)$ and $\GL_{k}(\mathbb C)$, respectively. Denote by $\mathbf E\to \mathcal M_0$ the holomorphic vector bundle that is determined by  the localy free sheaf $\mathcal E = \mathcal J/\mathcal J^2$, where $\mathcal J$ is the sheaf of ideals generated by odd elements in $\mathcal O_{\mathcal M}$. In
\cite[Proposition 5.2]{Oni} it was proven that $\mathbf E$ is a homogeneous vector bundle corresponding to the representation $\ph$ of $H$ such that
$$
\ph|R = \rho_{1}^*\otimes\rho_{2}.
$$
Since $(\tilde{\mathcal O}_{\mathcal M})_p\simeq\bigwedge^p\mathcal E$, we have to find the vector space of global sections of
$\bigwedge^p\mathbf E$. This vector bundle corresponds to the representation $\bigwedge^p\ph=
\bigwedge^p(\rho^*_{1}\otimes\rho_{2})$. We need to find dominat weights of $\bigwedge^p\ph$.

Let us choose the following  Cartan subalgebra 
$$
\mathfrak
t=\{\operatorname{diag}(\mu_1,\dots,\mu_n)\}
$$
in $\mathfrak g =
\mathfrak{gl}_n(\mathbb C)= \Lie (G)$.
For $p>0$ any weight of this representation has the following form:
$$
\Lambda=-\mu_{i_1}-\dots-\mu_{i_p}+\mu_{j_1}+\dots+\mu_{j_p},
$$
where 
$$
1\le i_1,\dots,i_p\le n-k\quad \text{and}\quad n-k+1\le j_1,\dots,j_p\le n.
$$
For $p=0$ the highest weight is equal to $0$. The weight $\Lambda=0$ is clearly dominant. For $p>0$ assume that $\mu^1 =
\mu_{i_a}$, where $i_a=\operatorname{max}\{i_1, \dots i_p\}$, and $\mu^2
= \mu_{j_b}$, where $j_b= \operatorname{min}\{j_1, \dots, j_p\}$. Then, $(\Lambda, \mu^1 - \mu^2)<0$. Therefore, the weight $\Lambda $ is not dominant and we have 
$$
(\tilde{\mathcal O}_{\mathcal M})_0(\mathcal M_0) = \mathbb C\quad \text{and}\quad (\tilde{\mathcal O}_{\mathcal M})_p (\mathcal M_0) =
\{0\} \quad \text{for} \quad p > 0.
$$

Now let us compute the space of global holomorphic functions  $\mathcal O_{\mathcal M}(\mathcal M_0)$. Clearly, $\mathcal
J^p(\mathcal M_0)=\{0\}$ for large $p$. For $p\ge 0$ we have an exact sequence
$$
0\to \mathcal J^{p+1}(\mathcal M_0)\longrightarrow \mathcal J^p(\mathcal M_0)\longrightarrow(\tilde{\mathcal
O}_{\mathcal M})_p(\mathcal M_0).
$$
By induction we see that $\mathcal J^p(\mathcal M_0)=\{0\}$ for $p>0$. Hence, for $p=0$ our exact sequence has the form:
$$
0\to \mathcal J^0(\mathcal M_0) = \mathcal O_{\mathcal M}(\mathcal M_0) \longrightarrow \mathbb C.
$$
Note that on any supermanifold we have constant functions.
Hence, $\mathcal O_{\mathcal M}(\mathcal M_0) =
\mathcal J^0(\mathcal M_0)=\mathbb C$.
Using Lemma 1 and induction, we get the result.$\Box$
\medskip

\section{Vector fields on $\Pi$-symmetric flag supermanifolds}

The Lie superalgebra of holomorphic vector fields on $\Pi$-symmetric super-Grassmannian $\mathbf{\Pi
Gr}_{n|n,k|k}$ was computed in \cite{Oni}. 

\medskip
\t\label{teor super=Grassmannians} {\sl Let $\mathcal M = \mathbf{\Pi Gr}_{n|n,k|k}$ and $(n,k)\ne (2,1)$. Then 
	$$
	\mathfrak v(\mathcal M)\simeq \mathfrak{q}_{n}(\mathbb C)/\langle E_{2n}\rangle. $$ 

Assume that $\mathcal M = \mathbf{\Pi Gr}_{2|2,1|1}$. Then 
$$
\mathfrak v(\mathcal M)\simeq
 \mathfrak g +\!\!\!\!\!\!\supset\langle z\rangle.
 $$ 
 Here $\mathfrak g = \mathfrak g_{-1}\oplus \mathfrak g_0\oplus \mathfrak g_{1}$ is a $\mathbb Z$-graded Lie superalgebra defined in the following way.
 $$
 \mathfrak g_{-1}= V, \quad \mathfrak g_{0}= \mathfrak{sl}_2(\mathbb C), \quad  \mathfrak g_{1}= \langle \d \rangle,
 $$
  where $V= \mathfrak{sl}_2(\mathbb C)$ is the adjoint $\mathfrak{sl}_2(\mathbb C)$-module, $[\mathfrak g_{0}, \mathfrak g_{1}] = \{0\}$ and $[\d , -]$ maps identicaly $\mathfrak g_{-1}$ to $\mathfrak g_{0}$.  Here $z$ is the grading operator of the $\mathbb Z$-graded Lie superalgebra $\mathfrak{g}$.

}
\medskip

\noindent{\it Proof.} The first part of the proof was given in \cite[Theorem 5.2]{Oni}. The second part follows from \cite[Theorem 4.2]{Oni}. 
The proof is complete.$\Box$

\medskip

\noindent {\bf Remark.} 
For completeness we describe the more general Onishchik result in our particular case $\mathcal M = \mathbf{\Pi Gr}_{2|2,1|1}$.
In local chart on $\mathbf{\Pi Gr}_{2|2,1|1}$
\begin{equation}\label{eq local chart GR_2,1}
\left(
\begin{array}{cc}
x&\xi\\
1&0\\
\xi&x\\0&1\end{array} \right)
\end{equation}
the Lie superalgebra of vector fields $\mathfrak v(\mathbf{\Pi Gr}_{2|2,1|1})$ has the following form:
\begin{align*}
\mathfrak g_{-1} = \langle \frac{\partial}{\partial \xi}, x\frac{\partial}{\partial \xi}, x^2\frac{\partial}{\partial \xi} \rangle, \quad  
&\mathfrak g_{0} = \langle \frac{\partial}{\partial x}, x^2\frac{\partial}{\partial x} + 2x\xi\frac{\partial}{\partial \xi}, x\frac{\partial}{\partial x} + \xi\frac{\partial}{\partial \xi} \rangle,\\
\d = \xi \frac{\partial}{\partial x},&\quad z=\xi \frac{\partial}{\partial \xi}.
\end{align*}

We will need another description of the Lie superalgebra of vector fields on $\mathbf{\Pi Gr}_{2|2,1|1}$. We have 
$$
\mathfrak v(\mathbf{\Pi Gr}_{2|2,1|1}) \simeq \mathfrak{q}_{2}(\mathbb C)/\langle E_{4}\rangle \oplus \langle z\rangle 
$$
as $\mathfrak g_{0}$-modules.
More precisely, in the local chart (\ref{eq local chart GR_2,1}) the isomorphism is given by
\begin{align*}
(\mathfrak{q}_{2}(\mathbb C)/\langle E_{4}\rangle)_{\bar 0}\simeq \mathfrak g_{0}& = \langle \frac{\partial}{\partial x}, x^2\frac{\partial}{\partial x} + 2x\xi\frac{\partial}{\partial \xi}, x\frac{\partial}{\partial x} + \xi\frac{\partial}{\partial \xi} \rangle, \\
(\mathfrak{q}_{2}(\mathbb C)/\langle E_{4}\rangle)_{\bar 1}&\simeq \langle \frac{\partial}{\partial \xi}, x\frac{\partial}{\partial \xi} + \xi\frac{\partial}{\partial x}, x^2\frac{\partial}{\partial \xi}, \xi\frac{\partial}{\partial x} \rangle, \\
z&=\xi \frac{\partial}{\partial \xi}.
\end{align*}

\medskip

From now on we use the following notations: 
$$
\mathcal M=\mathbf{\Pi
F}^{n|n}_{k|k}, \quad \mathcal B= \mathbf{\Pi
Gr}_{n|n,k_1|k_1}, \quad \mathcal S=\mathbf{\Pi
F}_{k'|k'}^{k_1|k_1},
$$
where $k'= (k_2,\ldots, k_r)$. We also assume that $r>1$.
 By Proposition \ref{prop Bash} and Theorem \ref{teor constant functions} the projection of the superbundle
$\mathcal M\to \mathcal B$ determines the homomorphism of Lie superalgebras 
$$
\mathcal{P}:\mathfrak v(\mathcal M)\to\mathfrak
v(\mathcal B).
$$
This projection is equivariant. Hence for the natural Lie superalgebra homomorphisms $\mu: \mathfrak{q}_{n}(\mathbb C)\to\mathfrak
v(\mathcal M)$ and $\mu_{\mathcal B}: \mathfrak{q}_{n}(\mathbb C)\to\mathfrak v(\mathcal B)$ we have 
$$
\mu_{\mathcal B} = \mathcal{P}\circ\mu.
$$
Note that for $r>1$ the base $\mathcal B$ cannot be isomorphic to $\mathbf{\Pi
	Gr}_{2|2,1|1}$. Therefore, by Theorem \ref{teor super=Grassmannians}, the homomorphisms $\mu_{\mathcal B}$ and hence the homomorphism
$\mathcal{P}$ is surjective. We will prove that $\mathcal{P}$ is injective. Hence, 
$$
\mu = \mathcal{P}^{-1}\circ\mu_{\mathcal B}
$$
is surjective and  
$$
\mathfrak v(\mathcal M)\simeq\mathfrak{q}_{n}(\mathbb C)/ \langle
E_{2n}\rangle .
$$

Let us study $\Ker\mathcal{P} \subset \mathfrak v(\mathcal M)$. In previous sections 
we constructed a localy free sheaf $\tilde{\mathcal
W}$ on $\mathcal B_0$. We have a natural action of $G =
\GL_n(\mathbb C)$ on the sheaf $\mathcal W$ that preserves the filtration (\ref{eq filtrations in W}) and induces the action 
on the sheaf $\tilde{\mathcal W}$. Hence, the vector bundle $\mathbf W_0\to \mathcal B_0$ corresponding to
$\tilde{\mathcal W}_0$ is homogeneous.
We use notations from the proof of Theorem \ref{teor constant functions}. Let us compute the representation of $H\subset G$ in the fiber of $\mathbf W_0$ 
over the point $o = H\in \mathcal B_0$. We will identify $(\mathbf W_0)_o$ with the Lie superalgebra of vector fields $\mathfrak v(\mathcal S)$
on $\mathcal S$.

Consider a local chart that contains $o$ on the
$\Pi$-symmetric super-Grassmannian $\mathcal B$. For example we can take the chart corresponding to $I_{1\bar 0} = \{n-k_1+1,\ldots,n\}$.
The coordinate matrix (\ref{eq coordinats matrices P-sym}) in this case has the following form
\begin{equation}\label{eq loc chart on B}
Z_{I_1} =\left(
\begin{array}{cc}
X_1&\Xi_1\\
E_{k_1}&0\\
\Xi_1&X_1\\0&E_{k_1}\end{array} \right). 
\end{equation}
Let us choose an atlas of $\mathcal M$ in a neighborhood of $o$  defined by certain $I_{s\bar 0},\; s = 2,\ldots,r$, see
(\ref{eq coordinats matrices P-sym}). In notations (\ref{eq_stab Grassmannian}) and (\ref{eq loc chart on B}) the group $H$ acts on
$Z_{I_1}$ in the following way:
$$
\left(
\begin{array}{cccc} A&0&0 &0\\C&B&0&0\\0&0&A&0\\0&0&C&B\end{array} \right)
Z_{I_1} = \left(
\begin{array}{cc} A X_1&A\Xi_1\\ C X_1 + B & C\Xi_1\\
A\Xi_1 & A X_1\\ C\Xi_1 & C X_1 + B \end{array} \right).
$$
Hence, for $Z_{I_2}$ we have
\begin{equation}\label{eq action over o}
\begin{split} \left(
\begin{array}{cc} C X_1 + B & C\Xi_1\\ C\Xi_1 & C X_1 + B\end{array} \right)&
\left(
\begin{array}{cc}  X_2 & \Xi_2 \\ \Xi_2 & X_2\end{array} \right) =\\
= &\left(
\begin{array}{cc} BX_2 + C X_1X_2 + C\Xi_1\Xi_2 &
B\Xi_2 + C X_1\Xi_2 + C\Xi_1 X_2\\
B\Xi_2 + C X_1 \Xi_2 + C \Xi_1 X_2 & B X_2 + C X_1X_2 + C\Xi_1\Xi_2
\end{array} \right).
\end{split}
\end{equation}

Note that the local coordinates  $Z_{I_s},\; s \ge 2$, determine the local coordinate on $\mathcal S$. To obtain the action of $H$ in the fiber $(\bf W_0)_o$ in these coordinates we put  $X_1 = 0, \; \Xi_1=0$ in (\ref{eq action over o}) and  modify
$Z_{I_s},\; s \ge 3$, accordingly.
We see that the nilradical of $H$ and the subgroup
$\GL_{n-k_1}(\mathbb C)$ of $R$ act trivially on
$\mathcal S$ and that  the subgroup $\GL_{k_1}(\mathbb C)\subset
R$ acts in the natural way. This means that $H$ acts as the even part of the Lie supergroup 
$Q_{k_1}(\mathbb C)$ on $\Pi$-symmetric flag supermanifold $\mathcal S$, see (\ref{eq action of Q}). 

Furthermore, by induction we assume that 
$$
\mathfrak v(\mathcal S)\simeq \mathfrak{q}_{k_1}(\mathbb C)/ \langle
E_{2n}\rangle\quad \text{or} \quad \mathfrak v(\mathcal S)\simeq \mathfrak{q}_{2}(\mathbb C) / \langle
E_{4}\rangle +\!\!\!\!\!\!\supset\langle z\rangle.
$$
 Then the induced action of
$\GL_{k_1}(\mathbb C)$ on $\mathfrak v(\mathcal S)$ coinsides with the adjoint action of the even part of  $Q_{k_1}(\mathbb
C)$. Standard computations lead to the following lemma, where we denote by $\Ad_{k_1}$ the adjoint representation of $\GL_{k_1}(\mathbb C)$ on $\mathfrak {sl}_{k_1}(\mathbb C)$
and  $1$ is the one dimensional trivial representation of $\GL_{k_1}(\mathbb C)$.

\medskip

\lem\label{lem representation of Gl_k_1 in v(S)} 
 {\sl The representation $\psi$ of $H$ in the fiber
$(\mathbf W_0)_o = \mathfrak v(\mathcal S)$ is completely reducible. If
$\mathfrak v(\mathcal S)\simeq \mathfrak{q}_{k_1}(\mathbb C) / \langle
E_{2n}\rangle $, 
then
\begin{equation}\label{eq rep of H in v(S) general}
 \psi|_{\mathfrak v(\mathcal S)_{\bar 0}} = \Ad_{k_1},\,\,\,
\psi|_{\mathfrak v(\mathcal S)_{\bar 1}} = \Ad_{k_1} + 1.
\end{equation}
If $\mathfrak v(\mathcal S)\simeq \mathfrak{q}_{2}(\mathbb C) / \langle
E_{4}\rangle +\!\!\!\!\!\!\supset\langle z\rangle$, 
then
\begin{equation}\label{eq rep of H in v(S) case (2,1)}
\psi|_{\mathfrak v(\mathcal S)_{\bar 0}} = \Ad_{k_1}+1,\,\,\,
\psi|_{\mathfrak v(\mathcal S)_{\bar 1}} = \Ad_{k_1} + 1.
\end{equation}

}

\medskip

Further we will use the chart on $\mathbf{\Pi F}_{k|k}^{n|n}$
defined by 
$I_{s\bar 0}$, where $I_{1\bar 0}$ is as above, and  
$$
I_{s\bar 0}=\{k_{s-1}-k_s+1,\ldots,k_{s-1}\}
$$
for $s\ge 2$.
The coordinate matrix of this chart have the following form
$$
Z_{I_s} = \left(
\begin{array}{cc} X_s&\Xi_s\\ E_{k_s}&0\\
\Xi_s&X_s\\0&E_{k_s}\end{array}\right),
$$
where again the local coordinate are $X_s=(x^s_{ij})$ and $\Xi_s=(\xi^s_{ij})$. We denote this chart by $\mathcal U$.

\medskip
\lem \label{lem some fundamental vector fields} {\sl The following vector fields in $\mathcal U$ 
$$
\frac{\partial}{\partial x^1_{ij}},\,\,\frac{\partial}{\partial
\xi^1_{ij}},\,\, u_{ij}+ \frac{\partial}{\partial x^2_{ij}},\,\,
v_{ij}+ \frac{\partial}{\partial \xi^2_{ij}}
$$
are fundamental. This is they are induced by the natural action of  $\Q_{n}(\mathbb C)$ on $\mathcal M$.
Here $u_{ij}$ and $v_{ij}$ are vector field that depend only on coordinates from $Z_{I_1}$. }
\medskip

\noindent{\it Proof.} Let us prove this statement for example for the vector field
$\frac{\partial}{\partial x^1_{11}}$. This vector field corresponds to the one-parameter subgroup $\operatorname{exp}(tE_{1,n-k_1+1})$.
Indeed, the action of this subgroup is given by
$$
\left(\begin{array}{cc} X_1&\Xi_1\\ E_{k_1}&0\\
\Xi_1&X_1\\0&E_{k_1}\end{array}\right)\mapsto
\left(\begin{array}{cc}\tilde X_1&\Xi_1\\ E_{k_1}&0\\
\Xi_1&\tilde X_1\\0&E_{k_1}\end{array}\right)\quad \text{and}\quad Z_{I_s}\mapsto Z_{I_s},\;
s\ge 2,
$$
where
$$
\tilde X_1=\left(
             \begin{array}{ccc}
               t+x_{11}^1 & \ldots & x_{1k_1}^1 \\
               \vdots & \ddots & \vdots \\
               x_{n-k_1,1}^1 & \ldots & x_{n-k_1, k_1}^1 \\
             \end{array}\right). \Box
$$

\medskip

Let us choose a basis $(v_q)$ of $\mathfrak v(\mathcal
S)$. In \cite{Vi} it was proven that any holomorphic vertical vector field on $\mathcal M$ can be written uniquely in the form 
\begin{equation}\label{eq form of vertical vector fields}
w=\sum_q f_qv_q, 
\end{equation}
where $f_q$ are holomorphic functions on $\mathcal U$ depending only on coordinates from $Z_{I_1}$. Further, we will need the following lemma:

\medskip
\lem \label{lem assume ker Pi ne 0} {\sl Assume that $\operatorname{Ker}\mathcal{P}\ne \{0\}$. Then there exists a vector field $w\in \operatorname{Ker}\mathcal{P}\setminus\{0\}$, such that 
$w=\sum\limits_q f_qv_q,$ where $f_q$ are holomorphic functions depending only on even coordinates from $Z_{I_1}$.}

\medskip

\noindent{\it Proof.} Assume that in (\ref{eq form of vertical vector fields}) there is a non-trivial vector field $w$ such that a function $f_q$ depends for example on
 $\xi_{ij}^1$. Then $w = \xi_{ij}^1 w'+ w''$, where $w'$
and $w''$ are vertical vector fields and their coefficients (\ref{eq form of vertical vector fields}) do not depend on $\xi_{ij}^1$, and $w'\ne 0$. Using Lemma \ref{lem some fundamental vector fields} and the fact that $\operatorname{Ker}\mathcal{P}$ is an ideal in
$\mathfrak v(\mathcal M)$, we see that 
$$
w' =
[w, \frac{\partial}{\partial \xi^1_{ij}}]\in
\operatorname{Ker}\mathcal{P}.
$$
Hence, we can exclude  all odd coordinates $\xi_{ij}^1$. $\Box$

\medskip
\noindent {\bf Corollary.} We have $$
(\operatorname{Ker}\mathcal{P}\ne \{0\}) \Longrightarrow (\mathcal W_0(\mathcal B_0) \ne \{0\}).
$$

\medskip

We will need the following well-known construction for holomorphic homogeneous vector bundles.
Let us take any homogeneous vector bundle $\mathbf E$ over $G/H$ and $x_0=H$. Assume that  $v_{x_0}\in \mathbf E_{x_0}$
is an $H$-invariant. We can construct the $G$-invariant section of $\mathbf E$ corresponding to $v_{x_0}$ in the following way.
We set  
$$
v_{x_1}:=g \cdot  v_{x_0}\in \mathbf E_{x_1},\quad \text{ where}\quad x_1= g x_0.
$$
Clearly this definition does not depend on  $g\in G$ such that $x_1= g x_0$. Indeed, assume that $x_1= g_1x_0$ and $x_1=
g_2x_0$. Then $g_2=g_1 h$, where $h\in H$. Hence, $$
g_2(v_{x_0})=
(g_1h)(v_{x_0})=g_1( v_{x_0}).
$$

We use this construction to express locally the $G$-invariant section of $\mathbf W_0$ corresponding to an $H$-invariant in $(\mathbf W_0)_{x_0}$. Let us take the following element $g\in G$:
$$
g= \left(
\begin{array}{cc}
E_{n-k_1}&A\\
0&E_{k_1}
\end{array}
\right) \times 
\left(
\begin{array}{cc}
E_{n-k_1}&A\\
0&E_{k_1}
\end{array}
\right),
$$
where $A$ is any complex matrix.
Then $g$ acts on $\mathbf{\Pi F}^{n|n}_{k|k}$ in the following way:
\begin{equation}\label{eq action of g}
\begin{split}
	\left(
	\begin{array}{cccc}
		E_{n-k_1}&A&0&0\\
		0&E_{k_1}&0&0\\
		0&0&E_{n-k_1}&A\\
		0&0&0&E_{k_1}\\
	\end{array}
	\right)
	\left(
	\begin{array}{cc}
		X^1&\Xi^1\\
		E_{k_1}&0\\
		\Xi^1&X^1\\
		0&E_{k_1}\\
	\end{array}
	\right) = \left(
	\begin{array}{cc}
		X^1+A&\Xi^1\\
		E_{k_1}&0\\
		\Xi^1&X^1+A\\
		0&E_{k_1}\\
	\end{array}
	\right),\\
	 Z_{I_s}= Z_{I_s},\quad s>1.
\end{split} 
\end{equation}
We see that $x_0=H$ has coordinates
$X_1= \Xi_1=0$. Clearly, $\{x_1\, | \,\,x_1 = gx_0\}$ is an open set in $\mathcal B_0$. Moreover, element $g$ does not modify fiber coordinates. Therefore, the corresponding to an $H$-invariant $v_{x_0}\in (\mathbf W_0)_{x_0}$ section $v$ is the constant section $x_1 \mapsto v_{x_0}$ over the open set $\{x_1\, | \,\,x_1 = gx_0\}$.

\medskip
\t \label{teor Ker Pi = 0} {\sl Assume that $r>1$ and $\mathfrak v(\mathcal S)\simeq \mathfrak{q}_{k_1}(\mathbb C) / \langle E_{2k_1}\rangle$. Then
$\operatorname{Ker} \mathcal{P}=\{0\}$ and 
$$
\mathfrak v(\mathbf{\Pi F}^{n|n}_{k|k})\simeq \mathfrak{q}_{n}(\mathbb C) / \langle E_{2n}\rangle.
$$}

\noindent{\it Proof.}
 First let us compute the vector space of global sections of $\mathbf W_0$. As in Theorem \ref{teor constant functions}, we use the Borel-Weyl-Bott Theorem.
The representation $\psi$ of $H$ in $(\mathbf W_0)_o$ is described in Lemma \ref{lem representation of Gl_k_1 in v(S)}. From (\ref{eq rep of H in v(S) general}) it follows that the highest weights of $\psi$
have the form: 
$$
\mu_{n-k_1+1}-\mu_{n} \quad  (\times 2)\quad \text{and}\quad  0.
$$
The first highest weight is not dominant because by definition of 
$\Pi$-symmetric flag supermanifolds $k_1<n$. The second highest weight is clearly dominant.
Therefore, the vector space of global sections of $\mathbf W_0$
is the irreducible $G$-module with highest weight $0$. Therefore, $\mathcal W_0(\mathcal B_0)\simeq \mathbb C$.

Let $v_{o}$ be the $H$-invariant element form   $\mathfrak v(\mathcal
S)$. It is defined by the following one-paremeter subsupergroup
$$
\beta(\tau)=\left( \begin{array}{cc}
E_{k_1}&\tau E_{k_1}\\
\tau E_{k_1}& E_{k_1}\\
\end{array}\right)
$$
in the Lie supergroup $Q_{k_1}$.
In our chart we have
\begin{equation}\label{eq H-invariant}
v_{o}=2\sum_{ij}\xi^2_{ij}\frac{\partial}{\partial x^2_{ij}}+u,
\end{equation}
where $u$ is a vector field depending only on coordinates   from $Z_{I_s}$, $s\ge 3$. Above we have seen that the corresponding to $v_o$ global section is a constant section. Therefore, 
the unique global  section $v$ of $\mathbf W_0$ in our chart has also the form (\ref{eq H-invariant}).

Assume  that $\operatorname{Ker}\mathcal{P}\ne \{0\}$ and $w\in
\operatorname{Ker}\mathcal{P}/\{0\}$ is as in Lemma \ref{lem assume ker Pi ne 0}. Clearly the vector fields $w$ and $\alpha(w)$, where $\alpha : \mathcal W \to \mathcal W_0$, have the same form in our chart.  Therefore, $w=av$ for some $a\in \mathbb C$.
Furthermore, 
$$
\alpha([w,v_{ij}+\frac{\partial}{\partial \xi^2_{ij}}])=
2a\frac{\partial}{\partial x^2_{ij}}.
$$
The commutator of these vector fields is an even vector field and it belongs to
 $\operatorname{Ker}\mathcal{P}$. Hence, $$
 \alpha([w,v_{ij}+\frac{\partial}{\partial \xi^2_{ij}}])= 0,
 $$
 because $\mathbf W_0$ has no global even sections.
Therefore, $a=0$ and the proof is complete.$\Box$

\medskip

Now consider the case when the fiber of superbundle $\mathcal M$
is isomorphic to $\mathbf{\Pi Gr}_{2|2,1|1}$.

\medskip
\t \label{teor main result special case} {\sl Assume that $r=2$ and $\mathcal S=\mathbf{\Pi Gr}_{2|2,1|1}$. Then 
	$$\mathfrak v(\mathbf{\Pi F}^{n|n}_{2,1|2,1}) \simeq \mathfrak{q}_{n}(\mathbb C) / \langle E_{2n}\rangle.
	$$ }

\noindent{\it Proof.} As in Theorem \ref{teor Ker Pi = 0}, let us compute the space of global sections of $\mathbf W_0$. From (\ref{eq rep of H in v(S) case (2,1)}) it follows that
the highest weights of $\psi$ have the form: 
$$
\mu_{n-k_1+1}-\mu_{n} \quad
(\times 2) \quad \text{and}\quad  0\quad (\times 2).
$$ 
By the Borel-Weyl-Bott Theorem we get 
$$
\mathcal W_0(\mathcal B_0)_{\bar 0}\simeq \mathcal W_0(\mathcal B_0)_{\bar 1}\simeq \mathbb
C.
$$
 A basic section of $\mathcal W_0(\mathcal B_0)_{\bar 1}$ was obtained in Theorem \ref{teor Ker Pi = 0}. It has the form $v=2\xi^2_{11}\frac{\partial}{\partial
x^2_{11}} $ in our case. Furthermore, we can take $z=\xi^2_{11}\frac{\partial}{\partial \xi^2_{11}}$ as a basic element of $1$-dimensional vector subspace in $\mathfrak v(\mathcal S)_{\bar 0}$ corresponding to the trivial representation $1$. Again in our local chart the unique even section of $\mathcal W_0(\mathcal B_0)_{\bar 0}$ has locally the form
$s=\xi^2_{11}\frac{\partial}{\partial \xi^2_{11}}$.

Let us take the vector field $w$ as in Lemma \ref{lem assume ker Pi ne 0}. The vector fields $\alpha(w)$ and $w$ have the same form in our chart.
Hence, $w=av+bs$, where $a,b\in \mathbb C$. Furthermore,
$$
\aligned 
\alpha([w,v_{11}+\frac{\partial}{\partial \xi^2_{11}}])=
\alpha(2a\frac{\partial}{\partial
x^2_{11}}+ b\frac{\partial}{\partial \xi^2_{11}})=
2a\frac{\partial}{\partial x^2_{11}}+b\frac{\partial}{\partial
\xi^2_{11}}\in \mathcal W_0(\mathcal B_0)=\langle v, z\rangle.
\endaligned
$$
In other words, $2a\frac{\partial}{\partial x^2_{11}}+b\frac{\partial}{\partial
	\xi^2_{11}}$ must be a linear combination of $v$ and $z$. Hence, $a=b=0$, and the proof is compete.$\Box$

\medskip

By induction we get our main result:

\medskip
\t\label{teor main result} {\sl Assume that $r>1$. Then 
	$\mathfrak v(\mathbf{\Pi F}^{n|n}_{k|k}) \simeq \mathfrak{q}_{n}(\mathbb C) / \langle E_{2n}\rangle.$ }

\section{The automorphism supergroup of $\mathbf{\Pi F}^{n|n}_{k|k}$}

Our main result has infinitesimal nature. However we can determine the connected component of the automorphism supergroup of $\mathbf{\Pi F}^{n|n}_{k|k}$. 
 Let us discuss this statement in details. 

Let us take a complex-analytic supermanifold $\mathcal M$ with a compact underlying space $\mathcal M_0$. Then the Lie superalgebra of holomorphic vector fields $\mathfrak v (\mathcal M)$ is finite dimensional. Denote by $\mathcal{A}ut (\mathcal M)_{\bar 0}$ the Lie group of even global automorphisms of $\mathcal M$. (The fact that $\mathcal{A}ut (\mathcal M)_{\bar 0}$ is a complex-analytic Lie group with the Lie algebra $\mathfrak v (\mathcal M)_{\bar 0}$ was proven in \cite{Bergner}.) Moreover, we have a natural holomorphic action of $\mathcal{A}ut (\mathcal M)_{\bar 0}$ on $\mathcal M$ (see \cite{Bergner}) and hence on $\mathfrak v (\mathcal M)$. Therefore, the pair $(\mathcal{A}ut (\mathcal M)_{\bar 0},\mathfrak v (\mathcal M))$ is a super Harish-Chandra pair. (See cite{ViLieSupergroup} for the definition of a super Harish-Chandra pair.) Using the equivalence of complex super Harish-Chandra pairs and complex Lie supergroups obtained in \cite{ViLieSupergroup} we determine the complex Lie supergroup $\mathcal{A}ut (\mathcal M)$. We call this Lie supergroup the {\it automorphism supergroup} of $\mathcal M$.  

Consider the case $\mathcal M = \mathbf{\Pi F}^{n|n}_{k|k}$.  Above we described a holomorphic action of $\GL_{n}(\mathbb C)=\Q_{n}(\mathbb C)_{\bar 0}$ on $\mathcal M$. In other words we have a homomorphism of Lie groups 
\begin{equation}\label{eq_homomorphism}
\Q_{n}(\mathbb C)_{\bar 0}\to \mathcal{A}ut (\mathcal M)_{\bar 0}.
\end{equation}
This homomorphism induces (almost always, see Theorems \ref{teor super=Grassmannians} and \ref{teor main result}) the isomorphism of Lie algebras $\mathfrak q_n(\mathbb C)_{\bar 0}/\langle E_{2n}\rangle$ and $\mathfrak v (\mathcal M)_{\bar 0}$, see Theorems \ref{teor super=Grassmannians} and \ref{teor main result}. In Section $2$ we have seen that the kernel of the homomorphism  (\ref{eq_homomorphism}) is equal to $\{\alpha E_{2n}\}$, where $\alpha \ne 0$, or to the center $\mathcal Z(\Q_{n}(\mathbb C)_{\bar 0})$ of $\Q_{n}(\mathbb C)_{\bar 0}$. Therefore, the connected component of the automorphism supergroup $\mathcal{A}ut^0 (\mathbf{\Pi F}^{n|n}_{k|k})$ is determined by the super Harish-Chandra pair 
$$
(\Q_{n}(\mathbb C)_{\bar 0}/ \mathcal Z(\Q_{n}(\mathbb C)_{\bar 0})
, \mathfrak q_n(\mathbb C)/\langle E_{2n}\rangle).
$$
In other words, 
$$
\mathcal{A}ut^0 (\mathbf{\Pi F}^{n|n}_{k|k})\simeq \Q_{n}(\mathbb C)/ \mathcal Z(\Q_{n}(\mathbb C)).
$$

In case $\mathcal M = \mathbf{\Pi Gr}_{2|2,1|1}$, the connected component of the automorphism supergroup $\mathcal{A}ut^0 (\mathbf{\Pi Gr}_{2|2,1|1})$ is given by the following super Harish-Chandra pair:
$$
(\Q_{2}(\mathbb C)_{\bar 0}/ \mathcal Z(\Q_{2}(\mathbb C)_{\bar 0})\times \mathbb C^*, \mathfrak v (\mathbf{\Pi Gr}_{2|2,1|1})),
$$
see Theorem \ref{teor super=Grassmannians} for a description of $\mathfrak v (\mathbf{\Pi Gr}_{2|2,1|1})$.

\noindent{\it Elizaveta Vishnyakova}

\noindent {Max Planck Institute for Mathematics, Bonn}

\noindent{\emph{E-mail address:}
	\verb"VishnyakovaE@googlemail.com"}


\begin{thebibliography}{99}
\bibitem[A]{ADima}  {\it  Akhiezer D. N.} Homogeneous complex manifolds. (Russian) Current problems in mathematics. Fundamental directions, Vol. 10 (Russian), 223-275, 283.

\bibitem[B]{Bash} {\it Bashkin M.A.} Contemporary problems in mathematics and informatics, vol. 3, Yaroslavl' State Univ., Yaroslavl' 2000, pp. 11–16.

\bibitem[BL]{BL} {\it Berezin F.A., Leites D.A.} Supermanifolds. Soviet
Math. Dokl. 16, 1975, 1218-1222.

\bibitem[BK]{Bergner} {\it Bergner H., Kalus M.} Automorphism groups of compact complex supermanifolds. arXiv:1506.01295

\bibitem[Kac]{Kac} {\it Kac V. G.} Classification of simple Lie superalgebras, Funktsional. Anal. i Prilozhen., 9:3 (1975), 91-92.

\bibitem[Man]{Man} {\it Manin Yu.I.} Gauge field theory and complex geometry, Grundlehren Math.
Wiss., vol. 289, Springer-Verlag, Berlin 1988, 1997.
	
	
\bibitem[Oni]{Oni} {\it Onishchik A.L.} Non-split supermanifolds associated with the cotangent bundle. Univer\-sit\'{e} de Poitiers,
D\'{e}partement de Math., N 109. Poitiers, 1997.


\bibitem[V1]{Vi} {Vishnyakova E.G.} Vector fields on flag supermanifolds (in Russian).
Sovremennye problemy mathematiki i informatiki, V. 8, Yaroslavl' State Univ., 2006, 11-23. 
	
\bibitem[V2]{ViPi-sym} {Vishnyakova E.G.} Vector fields on $\Pi$-symmetric flag supermanifolds (in Russian).
Vestnik TvGU, V. 7, Tver, 2007, 117-127. 

\bibitem[V3]{Viholom} {Vishnyakova E.G.} On holomorphic functions on a compact complex homogeneous supermanifold.
Journal of Algebra, Volume 350, Issue 1, 15 January 2012, Pages 174-196.

\bibitem[V4]{Vivector}{Vishnyakova E.G.} Lie superalgebras of vector fields on flag supermanifolds,  Russ. Math. Surv. 63 394, 2008.
	
\bibitem[V5]{ViLieSupergroup}{Vishnyakova E.G.}	 On complex Lie supergroups and split homogeneous supermanifolds. Transformation Groups 16 (2011), no. 1, 265-285.
	
\end{thebibliography}
\end{document}